\documentclass[mtplusscr,seceqn,dvips]{arxbj}
\usepackage{mathbh}
\usepackage{cursive}
\usepackage{euscript}
\usepackage{mathrsfs}

\aid{0}
\volume{13}
\issue{4}
\pubyear{2007}
\firstpage{1179}
\lastpage{1194}
\doi{10.3150/07-BEJ6180}

\makeatletter
\newtheorem{theorem}{Theorem}[section]
\newtheorem{corollary}{Corollary}[section]
\newproclaim{definition}{Definition}[section]
\newremark{remark}{Remark}[section]
\newtheorem{lemma}{Lemma}[section]
\newremark{warning}{Warning}
\makeatother

\begin{document}
\begin{frontmatter}

\title{The delta method for analytic functions of random operators
with application to functional data}
\runtitle{The delta method for analytic functions}

\begin{aug}
\author[1]{\fnms{J.} \snm{Cupidon}\thanksref{1}},
\author[1]{\fnms{D.S.} \snm{Gilliam}\corref{}\thanksref{1}\ead[label=e1]{gilliam@math.ttu.edu}},
\author[2]{\fnms{R.} \snm{Eubank}\thanksref{2}}
\and
\author[1]{\fnms{F.} \snm{Ruymgaart}\thanksref{1}}
\runauthor{J. Cupidon, D.S. Gilliam, R. Eubank and F. Ruymgaart}
\address[1]{Mathematics and Statistics,
Texas Tech University,
Lubbock, TX 79409, USA.\\
\printead{e1}}
\address[2]{Mathematics,
Arizona State University,
Tempe, AZ, 85278-1804, USA.}
\end{aug}

\received{\smonth{11} \syear{2006}}
\revised{\smonth{6} \syear{2007}}

%
\begin{abstract}
In this paper, the asymptotic distributions of estimators for the
regularized functional
canonical correlation and variates of the population are derived. The
method is based on the possibility
of expressing these regularized quantities as the maximum eigenvalue
and the corresponding
eigenfunctions of an associated pair of regularized operators, similar
to the Euclidean case.
The known weak convergence of the sample covariance operator, coupled
with a delta-method for
analytic functions of covariance operators, yields the weak convergence
of the pair of
associated operators. From the latter weak convergence, the limiting
distributions of the
canonical quantities of interest can be derived with the help of some
further perturbation theory.
\end{abstract}

%

\end{frontmatter}
%

\section{Introduction }\label{sec1}

This paper deals with the asymptotic distribution theory of functional
canonical correlations and their variates. Although tailored to these
particular problems, the methodology is of a generic character and may
also apply to questions regarding the asymptotic distribution of other
statistics used in functional data analysis. The problem will
be formulated in a general Hilbert space setting where the Hilbert
space is tacitly assumed to be infinite-dimensional and separable.

In this infinite-dimensional case, some difficulties regarding the
definition of the sample canonical correlation have already been
observed in Leurgans \textit{et al.} (\citeyear{lemosi93}). The authors of that paper argue that some
kind of smoothing or regularization is indispensable when dealing with
the sample canonical correlation. These difficulties are
essentially due to the fact that the sample covariance operator has a
so-called finite-dimensional kernel (Riesz and Sz.-Nagy (\citeyear{rina90})), while
acting on an infinite-dimensional space. Leurgans \textit{et al.} (\citeyear{lemosi93})
realize smoothing by introducing a roughness penalty term. Although there is
a connection between Tikhonov regularization of inverse operators
(employed in this paper) and the use of penalty terms, the relation
with the roughness penalty cannot be established within the present context
of our paper. He \textit{et al.} (\citeyear{hemuwa04}) apply dimension
reduction/augmentation at the level of the actual data and base the empirical
canonical correlation on these modified data. This approach differs
considerably from ours, which is based on regularization of the
canonical correlation itself. The results in He \textit{et al.} (\citeyear{hemuwa04}) are
for fixed sample size and the asymptotics in Leurgans \textit{et al.}
(\citeyear{lemosi93}) remain restricted to consistency.

In \citeyear{cueugiru06} it has been observed that the population
canonical correlation, although well defined in principle, is, in
general, a supremum of a certain functional, rather than a maximum, so
that a maximizer (i.e., a pair of canonical variates) may not always
exist in the ambient Hilbert space. Another deficiency is that, even if
the canonical correlation corresponds to a maximum and canonical
variates do exist, these quantities cannot be interpreted as
the maximum eigenvalue and corresponding eigenvector of a pair of
associated operators, as is true in the Euclidean setting. The
development in \citeyear{cueugiru06} shows that all of these
deficiencies of the population canonical correlation can be remedied if
a modification is employed, based on regularization of the inverses of
the operators involved. Also, some relations between the actual
population quantities and their regularized versions are established in
that paper.

The present approach to finding the asymptotic distribution of the
regularized sample canonical correlation and its variates hinges
to a great extent on the interpretation of both the regularized sample and the
regularized population quantities as spectral characteristics of
associated pairs of operators. In Section~\ref{sec4} of this paper, the
asymptotic distribution of a regularized version of the sample
canonical correlation and its variates will be derived. In the
Euclidean case, where regularization is not needed, this approach has
been pursued in Ruymgaart and Yang (\citeyear{ruya97}), exploiting certain results in
Watson (\citeyear{watson83}).

One of the main tools needed to derive the desired asymptotics is a
delta-method for analytic functions of certain random operators (more
specifically, sample covariance operators). This delta-method might be of
independent interest and is considered in Section~\ref{sec3}. It is
based on the existence of a Fr\'echet derivative of an analytic
function of a compact, strictly positive Hermitian operator,
tangentially to the space of all compact Hermitian operators. Because
we cannot make the simplifying assumption that the increments commute
with the operator at which the function is evaluated, the expression
for the Fr\'echet derivative requires an extra correction term. The
delta-method yields the asymptotic distribution of the associated
operators, from which the asymptotics of their eigenvalues and
eigenvectors can be derived in a similar manner as in Dauxois \textit{et al.} (\citeyear{daporo82}).

As has been observed above, without regularization, the population
canonical variates do not, in general, exist and, consequently, it seems
appropriate to maintain a fixed level of regularization for suitable
asymptotics. Mathematically, a fixed level of regularization leads to
root-sample-size asymptotics. When the regularization parameter tends
to zero, however, this rate will depend on the (typically unknown)
eigenvalues of the covariance operator.

In Section \ref{sec2}, some basic notation and definitions are introduced.
For practical implementation of the results of Section~\ref{sec4}, the
estimation of unknown parameters will be needed, an issue addressed in
Section~\ref{sec5}. An example and some further comments are given in
Section~\ref{sec6}. The mathematical results for perturbation of
compact, positive Hermitian operators that, in particular, yield the Fr\'
echet derivative are reviewed without proof in the \hyperref[secA]{Appendix}.

\section{Basic notation, definitions and assumptions}\label{sec2}

Let $(\Omega ,\EuScript{F},\mathbb{P})$ denote a probability space,
$\mathbb{H}$ an infinite dimensional, separable Hilbert space
with inner product $\langle {\cdot},{\cdot} \rangle $, norm $\|
\cdot\|$
and $\sigma $-field of Borel sets $\EuScript{B}_\mathbb{H}$, and let
$X\dvtx  \Omega \rightarrow \mathbb{H}$ be a random element in
$\mathbb{H}$, that is, an $(\EuScript{F},\EuScript{B}_\mathbb
{H})$-measurable mapping. Throughout, it will be required that
%
\begin{equation}
\label{e2.1}
\mathbb{E}\|X\|^4 <\infty .
\end{equation}
Under this condition, the mean $\mathbb{E}X = \mu\in\mathbb{H}$
exists, meaning that (Laha and Rohatgi (\citeyear{laro79}))
%
\begin{equation}
\label{e2.2}
\mathbb{E}\langle {f},{X} \rangle = \langle {f},{\mu} \rangle\qquad
 \forall f\in\mathbb{H}.
\end{equation}

Under assumption (\ref{e2.1}), the covariance operator $\Sigma$ of $X$
also exists. It is known to be uniquely determined by the relation

\begin{equation}
\label{e2.3}
\mathbb{E}\langle {f},{X-\mu} \rangle \langle {X-\mu},{g}
\rangle =\mathbb{E}\big\langle f,\big((X-\mu)\otimes(X-\mu)\big
)g \big\rangle = \langle {f},{\Sigma g} \rangle\qquad\forall f,g\in
\mathbb{H},
\end{equation}
where ``$\otimes$'' denotes the tensor product in $\mathbb{H}$. We will
also write
%
\begin{equation}
\label{e2.4}
\Sigma= \mathbb{E}(X-\mu)\otimes(X-\mu).
\end{equation}
Such a covariance operator is nonnegative Hermitian and has finite
trace $\mathbb{E}\|X\|^2$, so it is also compact. We will
therefore assume, without real loss of generality, that
%
\begin{equation}
\label{e2.5}
\Sigma\mbox{ is strictly positive, that is, }  \langle {f},{\Sigma f}
\rangle >0\qquad  \forall f\neq0
\end{equation}
and hence that $\Sigma$ is injective. It is well known that $\Sigma $
has spectral representation
%
\begin{equation}
\label{e2.6}
\Sigma = \sum_{k=1}^\infty \lambda _k P_k,
\end{equation}
where $\lambda _1>\lambda _2 > \cdots\downarrow0$ are the
eigenvalues of $\Sigma $ and $P_1, P_2, \ldots$ the projections onto
the corresponding finite-dimensional eigenspaces.

Let $\mathcal{L}$ denote the Banach space of all bounded linear
operators that map $\mathbb{H}$ into itself. The ordinary operator
norm in $\mathcal{L}$ will be denoted by $\|\cdot\|$ without confusion.
Of particular importance in this paper, however, is the subspace
$\mathcal{L}(\mbox{HS})$ of all Hilbert--Schmidt operators. This space
becomes a separable Hilbert space when it is endowed with the inner product
%
\begin{equation}
\label{e2.7}
\langle {U},{V} \rangle_{\mathrm{HS}}= \sum_{k=1}^\infty \langle
{Ue_k},{Ve_k} \rangle ,\qquad U, V \in\mathcal{L}(\mbox{HS}),
\end{equation}
where $e_1, e_2, \ldots$ is an orthonormal basis of $\mathbb{H}$.
This inner product does not depend on the choice of basis; see Lax
(\citeyear{lax00}). The norm and tensor product in $\mathcal{L}(\mbox{HS})$ will
be denoted by $\|\cdot\|_{\mathrm{HS}}$ and $\otimes_\mathrm{HS}$,
respectively.

The space $\mathcal{L}(\mbox{HS})$ is important for the study of weak
convergence of the sample covariance operator. At this point, let us
simply note that $\Sigma\in\mathcal{L}(\mbox{HS})$ and that $(X-\mu
)\otimes(X-\mu)$ is a random element in $\mathcal{L}(\mbox{HS})$.
As a random element in this Hilbert space, it has its own covariance
operator; this operator exists due to condition (\ref{e2.1}) and
can easily be seen to equal (cf. (\ref{e2.3}) and (\ref{e2.4}))
%
\begin{eqnarray}\label{e2.8}
&& \mathbb{E}\{ (X-\mu)\otimes(X-\mu)-\Sigma \}
\otimes_{\mathrm{HS}} \{ (X-\mu)\otimes(X-\mu)-\Sigma
\} \nonumber \\
&&\qquad = \mathbb{E}\{ (X-\mu)\otimes(X-\mu) \} \otimes
_{\mathrm{HS}} \{ (X-\mu)\otimes(X-\mu) \} -\Sigma
\otimes_{\mathrm{HS}}\Sigma \\
&&\qquad= \Sigma_\mathrm{HS}.\nonumber
\end{eqnarray}

Next, let us suppose that $\mathbb{H}_1$ and $\mathbb{H}_2$ are two
closed subspaces of $\mathbb{H}$ such that
%
\begin{equation}
\label{e2.9}
\mathbb{H}=\mathbb{H}_1\oplus\mathbb{H}_2,\qquad\mathbb{H}_1\perp
\mathbb{H}_2.
\end{equation}
Denote the orthogonal projection of $\mathbb{H}$ onto $\mathbb{H}_j$
by $\Pi_j$, let $X_j = \Pi_jX$, $\mu_j =\Pi_j \mu$ and let $\Sigma
_{jk}$ denote the restriction of $\Sigma$ to $\mathbb{H}_k$ and $\mathbb
{H}_j$, that is,
%
\begin{equation}
\label{e2.10}
\Sigma_{jk} = \Pi_j\Sigma\Pi_k,\qquad j,k = 1,2.
\end{equation}
Because the $\Pi_j$ are bounded and $\Sigma$ is Hilbert--Schmidt (and
hence compact), each operator $\Sigma_{jk}$ is still Hilbert--Schmidt
(and hence compact).
In addition, the $\Sigma_{jj}$ are strictly positive Hermitian. Let us
also note that
%
\begin{equation}
\label{e2.11}
\Sigma_{12}^*=(\Pi_1\Sigma\Pi_2)^* = \Pi_2 \Sigma\Pi_1 = \Sigma_{21}.
\end{equation}
Similarly to (\ref{e2.6}), $\Sigma_{jj}$ has a spectral representation
of the form
%
\begin{equation}
\label{e2.12}
\Sigma_{jj} = \sum_{k=1}^\infty \lambda _{jk}P_{jk},\qquad j=1,2,
\end{equation}
where $\lambda _{j1}>\lambda _{j2}> \ldots\downarrow0$ are the
eigenvalues of $\Sigma_{jj}$ and $P_{j1}, P_{j2}, \ldots$ the
projections onto the corresponding finite dimensional eigenspaces.

Suppose, now, that we are given a random sample $X_1,X_2, \ldots, X_n$
of independent copies of~$X$. The usual estimators of $\mu$ and
$\Sigma$ are
%
\begin{equation}
\label{e2.13}
\overline{X} = \frac{1}{n} \sum_{i=1}^n X_i,\qquad \widehat {\Sigma
} = \frac{1}{n} \sum_{i=1}^n (X_i-\overline{X})\otimes
(X_i-\overline{X}),
\end{equation}
respectively.
This operator $\widehat {\Sigma}$ has all of the properties of
$\Sigma$, including its being of Hilbert--Schmidt type, except that it has a so-called
finite-dimensional kernel (Riesz and Sz. Nagy (\citeyear{rina90})) with a
range of dimension at most $n$. Hence, this operator can never be
injective, not even when $\Sigma$ is (as we assume).
The fact that $\widehat {\Sigma}$ is not injective is the source of
difficulties associated with defining the sample principal canonical correlation
that turns out to always be $1$, as has been pointed out by Leurgans
\textit{et al.} (\citeyear{lemosi93}). These authors state that regularization is
indispensable in the sample case.

The canonical correlation concept considered here can also be viewed
from the perspective of Hilbert-space-indexed processes (e.g., Parzen
(\citeyear{pa70})) corresponding to $\mathbb{H}$ inner products involving the
random elements $X_j=\Pi_j X$, $j=1,2$. Thus, it has direct ties to
(functional) analysis of variance and discriminant analysis that
parallel the relationship between these methods for classical
multivariate analysis (e.g., Kshirsagar (\citeyear{ks72}), Eubank and Hsing (\citeyear{euhs06})
and Shin (\citeyear{shin06})). The necessity of regularization in this context
follows from results in Bickel and Levina (\citeyear{bile04}), while the use of
regularized discriminant analysis methods with functional data has been
explored by Hastie \textit{et al.} (\citeyear{habuti95}).


\citeyear{cueugiru06} argue that regularization is expedient, even
when the population canonical correlation is considered, because,
without it, canonical variates may not exist and the relation with the
spectral characteristics of an associated pair of operators is lost.
Hence, in this paper, both the sample and the population canonical
correlation will be regularized and compared at the same fixed, but
arbitrary, level of the regularization parameter.

In order to specify the regularization that will be employed here, let
us replace $\Sigma $ with $\alpha I+\Sigma $ and $\widehat
{\Sigma }$ with $\alpha I+\widehat {\Sigma }$, where $I$ is the
identity operator and $\alpha >0$. Let us also replace $\Sigma _{jk}$
and $\widehat {\Sigma }_{jk}$ with
%
\begin{eqnarray}
\Pi_j(\alpha I+\Sigma )\Pi_k =
\cases{ (\alpha I_j + \Sigma _{jj}), &$\quad j=k$,\cr \Sigma
_{jk},&$\quad j\neq k,$}\label{e2.14}\\
\Pi_j(\alpha I+\widehat {\Sigma })\Pi_k =
\cases{(\alpha I_j + \widehat {\Sigma }_{jj}), &$\quad j=k$,\cr
\widehat {\Sigma }_{jk},&$\quad j\neq k$,}\label{e2.15}
\end{eqnarray}
respectively, where $I_j =\Pi_j$ is essentially the identity operator
restricted to $\mathbb{H}_j$. Let us write $\mathbb{H}_1^0=\mathbb
{H}_1 \backslash\{0\}$, $\mathbb{H}^0_2=\mathbb{H}_2\backslash\{0\}
$ for brevity.
\begin{definition}\label{def2.1}
Fix $\alpha >0$. The regularized squared principal canonical
correlation (RSPCC) for the population is defined as
%
\begin{equation}
\label{e2.16}
\rho ^2=\rho ^2(\alpha ) = \mathop{\max_{f_1\in\mathbb{H}^0_1}}_
{f_2\in\mathbb{H}^0_2} \frac{\langle {f_1},{\Sigma _{12}f_2}
\rangle ^2}{\langle {f_{1}},{(\alpha I_1+\Sigma _{11})f_1} \rangle
\langle {f_2},{(\alpha I_2 +\Sigma _{22})f_2} \rangle } .
\end{equation}
Its sample analogue is $\widehat {\rho }^2 = \widehat {\rho
}^2(\alpha )$, obtained from (\ref{e2.16}) by replacing $\Sigma
_{jk}$ with $\widehat {\Sigma }_{jk}$. Pairs of maximizers will
be respectively denoted by $f_1^*=f_{1\alpha }^*$, $f_2^*=f_{2\alpha
}^*$ for the population and by $\widehat {f}_1=\widehat {f}_{1\alpha
}$, $\widehat {f}_2=\widehat {f}_{2\alpha }$ for the sample. The
corresponding canonical variates are
%
\begin{equation}
\label{e2.17}
\langle {X},{f^*_j} \rangle ,\  \langle {X},{\widehat {f}^*_j}
\rangle ,\qquad j=1,2.
\end{equation}
\end{definition}
\begin{warning*} Since, throughout the sequel, $\alpha >0$ will
be arbitrary, but fixed, the dependence on $\alpha $ is henceforth
suppressed in the notation.
\end{warning*}

Several properties have been shown in \citeyear{cueugiru06}, in
particular, that, for $\alpha >0$, a maximizer always exists. This can, in
fact, be seen as an implication of the following result of that paper.
Define the operators ($\alpha >0$)
%
\begin{eqnarray}
R_1 &=& (\alpha I_1 + \Sigma _{11})^{-1/2}\Sigma _{12}(\alpha I_2
+\Sigma _{22})^{-1} \Sigma _{21}(\alpha I_1+\Sigma _{11})^{-1/2},\label{e2.18}\\
R_2 &=& (\alpha I_2 + \Sigma _{22})^{-1/2}\Sigma _{21}(\alpha I_1
+\Sigma _{11})^{-1} \Sigma _{12}(\alpha I_2+\Sigma _{22})^{-1/2}\label{e2.19}
\end{eqnarray}
and their sample analogues $\widehat {R}_1$ and $\widehat {R}_2$.
Since all factors defining these operators are bounded, with $\Sigma
_{12}$ and $\Sigma _{21}$ or their sample analogues even
Hilbert--Schmidt (and hence compact) it follows that these operators are
also Hilbert--Schmidt (and hence compact). It will be assumed that
%
\begin{equation}
\label{e2.20}
\cases{R_{1}\mbox{ and }R_{2}\mbox{ have a largest eigenvalue }
\mbox{with one-dimensional eigenspace }\cr \mbox{generated by } f^*_{1}\mbox{ and }f^*_{2},
\mbox{respectively, where }\|f^*_{1}\|=\|f^*_{2}\|=1.}
\end{equation}
\begin{theorem}\label{thm2.1}
For $\alpha >0$, we have
%
\begin{equation}
\label{e2.21}
\rho ^2 = \mbox{ largest eigenvalue of } R_j = \langle
{f^*_j},{R_jf^*_j} \rangle
\end{equation}
for $j=1,2$. A similar result holds true for $\widehat {\rho }^2$.
\end{theorem}

The maximizers or canonical variates are essentially unique if the
eigenspaces corresponding to this maximal eigenvalue are
one-dimensional. The same properties hold true for the sample analogue.

\section{A delta-method for analytic functions of the sample
covariance operator}\label{sec3}

Assuming (\ref{e2.1}), Dauxois \textit{et al.} (\citeyear{daporo82}) have shown the
fundamental result
%
\begin{equation}
\label{e3.1}
\sqrt{n}(\widehat {\Sigma }-\Sigma )\mathop{\stackrel{d}{\rightarrow}}
\mathcal{G},\qquad \mbox{as }  n\rightarrow \infty , \mbox{
in } \mathcal{L}(\mbox{HS}),
\end{equation}
where $\mathcal{G}$ is a zero-mean Gaussian random element in the
Hilbert space $\mathcal{L}(\mbox{HS})$ with covariance operator
%
\begin{equation}
\label{e3.2}
\mathbb{E}\, \mathcal{G}\otimes_{\mathrm{HS}} \mathcal{G}= \Sigma
_{\mathrm{HS}},
\end{equation}
as defined in (\ref{e2.8}). The continuous mapping theorem immediately yields
 that
%
\begin{equation}
\label{e3.3}
\sqrt{n}(\widehat {\Sigma }_{jk}-\Sigma _{jk})\mathop{\stackrel{d}{\rightarrow}}
\Pi_j\mathcal{G}\Pi_k=\mathcal{G}_{jk} ,\qquad \mbox{as }  n\rightarrow \infty , \mbox{ in } \mathcal{L}(\mbox{HS}).
\end{equation}

Let $D\subset\mathbb{C}$ be the open domain in the complex plane
defined by
%
\begin{equation}
\label{e3.4}
D=\biggl\{ z\in\mathbb{C}\dvtx \min_{0\leq x\leq\|\Sigma \|
}|z-x|< \tfrac{1}{2} \alpha \biggr\},
\end{equation}
where $\alpha >0$ is the regularization parameter. This domain can be
used for all the specific functions we need to consider. It seems
worthwhile, however, to first consider an arbitrary function
%
\begin{equation}
\label{e3.5}
\varphi\dvtx D \rightarrow \mathbb{C},\qquad\mbox{analytic on }  D.
\end{equation}
As in the \hyperref[secA]{Appendix}, let $\mathcal{C}_H$ denote the class of all compact
Hermitian operators on $\mathbb{H}$ and $\mathcal{L}_H$ the class of
all bounded Hermitian operators. Let us consider the operator $\varphi
(\Sigma +\EuScript{P})$ in $\mathcal{L}_H$, for $\EuScript{P}$ in
$\mathcal{C}_H$ with $\|\EuScript{P}\|<\frac{1}{3}\alpha $. This
operator-valued function has a Fr\'echet derivative at $\Sigma $,
tangentially to $\mathcal{C}_H$, denoted by $\varphi '_\Sigma $ and
given by (\ref{eA.6}).
This operator $\varphi '_\Sigma\dvtx \mathcal{C}_H\rightarrow
\mathcal{L}_H$ is bounded in the usual operator norm.

If $\mathcal{L}_H(\mbox{HS}) \subset\mathcal{C}_H$ is the subspace
of all Hermitian Hilbert--Schmidt operators, we even have
%
\begin{equation}
\label{e3.6}
\varphi _\Sigma'\dvtx \mathcal{L}_H(\mbox{HS}) \rightarrow
\mathcal{L}_H(\mbox{HS}),\qquad\mbox{bounded in } \|\cdot\|
_{\mathrm{HS}}.
\end{equation}
To see this, take $\EuScript{P}\in\mathcal{L}_H(\mbox{HS})$ and
observe that
%
\begin{equation}
\label{e3.7}
\|\varphi '_\Sigma \EuScript{P}\|^2_{\mathrm{HS}} = \sum
_{k=1}^\infty \|\varphi '_\Sigma \EuScript{P}e_k\|^2
\leq\|\varphi '_\Sigma \|^2 \sum_{k=1}^\infty \|\EuScript{P}e_k\|^2
= \|\varphi '_\Sigma \|^2 \|\EuScript{P}\|^2_{\mathrm{HS}} <\infty ,
\end{equation}
exploiting the boundedness of $\varphi '_\Sigma $ in the usual
operator norm. It is well known (Lax (\citeyear{lax00})) that
%
\begin{equation}
\label{e3.8}
\|T\|\leq\|T\|_{\mathrm{HS}},\qquad T\in\mathcal{L}(\mbox{HS}).
\end{equation}
We are now ready to establish a ``delta-method'' for random operators.
For random matrices, the result follows from Watson (\citeyear{watson83}) and can be
found in Ruymgaart and Yang (\citeyear{ruya97}).
\begin{theorem}\label{thm3.1}
If (\ref{e2.1}) is satisfied, it then follows that
%
\begin{equation}
\label{e3.9}
\sqrt{n}\{ \varphi (\widehat {\Sigma })-\varphi (\Sigma
)\} \mathop{\stackrel{d}{\rightarrow}}\mathcal{H},\qquad \mbox{as }
n\rightarrow \infty , \mbox{ in } \mathcal{L}(\mathrm{HS}),
\end{equation}
where $\mathcal{H}$ is the zero-mean Gaussian random element of
$\mathcal{L}(\mathrm{HS})$ given by
%
\begin{equation}
\label{e3.10}
\mathcal{H} = \varphi'_\Sigma \mathcal{G}= \sum_{j\geq1}
\varphi '(\lambda _j) P_j \mathcal{G}P_j
+ \mathop{\sum\sum}_{j\neq k} \frac{\varphi (\lambda _k)-\varphi
(\lambda _j)}{\lambda _k-\lambda _j} P_j \mathcal{G}P_k,
\end{equation}
with $\mathcal{G}$ given in (\ref{e3.1}).
\end{theorem}
\begin{pf*}{Proof}
Let us consider $\widehat {\EuScript{P}} = \widehat {\Sigma
}-\Sigma $ as a random perturbation (cf. Dauxois \textit{et al.} (\citeyear{daporo82}),
Watson (\citeyear{watson83})) and note that, by (\ref{e3.1}) and (\ref{e3.8}), we have
$\|\widehat {\EuScript{P}}\|\leq\|\widehat {\EuScript{P}}\|_{\mathrm
{HS}} =\EuScript{O}_p(n^{-1/2})$ as $n\rightarrow \infty $.
This implies that, for numbers $n^{-1/2} \ll\epsilon _n \ll n^{-1/4}$
we have
%
\begin{equation}
\label{e3.11}
\mathbb{P}(\Omega _n) = \mathbb{P}\{\omega \in\Omega\dvtx \|
\widehat {\EuScript{P}}(\omega )\|< \epsilon _n\} \rightarrow 1\qquad\mbox{as }  n\rightarrow \infty .
\end{equation}
According to Theorem \ref{thmA.1} and (\ref{e3.12}), we have, for $n$
sufficiently large,
%
\begin{eqnarray}
\label{e3.12}
\sqrt{n}\{ \varphi (\widehat {\Sigma })-\varphi (\Sigma
)\} & =& \sqrt{n}\{ \varphi (\widehat {\Sigma
})-\varphi (\Sigma )\}\mathbf{1}_{\Omega _n} + \sqrt
{n}\{ \varphi (\widehat {\Sigma })-\varphi (\Sigma )\}
\mathbf{1}_{\Omega _n^c}\nonumber \\
&=& \sqrt{n} \{\varphi '_\Sigma \widehat {\EuScript{P}} +
\EuScript{O}(\|\widehat {\EuScript{P}}\|^2) \}\mathbf
{1}_{\Omega _n}+\mathscr{o}_p(1) \\
&=& \varphi '_\Sigma \bigl(\sqrt{n}(\widehat {\Sigma }-\Sigma )
\bigr)+\mathscr{o}_p(1).\nonumber
\end{eqnarray}
The results in the theorem follow from (\ref{e3.12}) by applying
(\ref{e3.1}) once more, in conjunction with (\ref{e3.6}) and the
continuous mapping theorem.
\end{pf*}
\begin{remark}\label{rem3.1}
The double sum in (\ref{e3.1}) is, in fact, a correction term that is
needed because we may not assume that the ``increments'' $\widehat
{\EuScript{P}}=\widehat {\Sigma } -\Sigma $ and $\Sigma $ commute;
see also Remark \ref{remA.1}.
\end{remark}

In order to obtain asymptotic distributions for functional canonical
correlations and variates, Theorem \ref{thm3.1} will be employed for
the specific functions
%
\begin{equation}
\label{e3.13}
\varphi _p(z) = (\alpha + z)^{-p/2},\qquad z\in D, \ \ p=1,2.
\end{equation}
These functions are indeed analytic on $D$. For brevity, let us simply
write $\varphi '_{p,j}$ for the Fr\'echet derivative evaluated at
$\Sigma _{jj}$. It is immediate from (\ref{e2.9}) that $\|\Sigma
_{jj}\|\leq\|\Sigma \|$ and therefore the domain $D$ can still be
used for $\Sigma _{jj}$. The following corollary is immediate from
these remarks, (\ref{e3.3}) and Theorem \ref{thm3.1}.
\begin{corollary}\label{cor3.1}
With $\varphi _p$ as in (\ref{e3.13}), we have, for $j=1,2$,
%
\begin{equation}
\label{e3.14}
\sqrt{n}\{\varphi _p (\widehat {\Sigma }_{jj}) -\varphi
_p(\Sigma _{jj}) \}\mathop{\stackrel{d}{\rightarrow}} \varphi'_{p,j}\mathcal{G}_{jj},
\end{equation}
where the limit is a zero-mean Gaussian random element in $\mathcal
{L}(\mathrm{HS})$ and, more explicitly,
\begin{eqnarray}
\label{e3.15}
\varphi '_{p,j} \mathcal{G}_{jj}  &=& -\frac{p}{2} \sum_{k\geq
1}\frac{1}{(\alpha +\lambda _{jk})^{(p+2)/2}} P_{jk}\mathcal
{G}_{jj} P_{jk}\nonumber\\[-8pt]\\[-8pt]
&&\hspace*{-5pt}{}+\mathop{\sum\sum}_{m\neq n} \frac{(\alpha +\lambda
_{jm})^{p/2}-(\alpha +\lambda _{jn})^{p/2}}{(\lambda _{jn}-\lambda
_{jm}) (\alpha +\lambda _{jm})^{p/2} (\alpha +\lambda _{jn})^{p/2}}
P_{jm}\mathcal{G}_{jj} P_{jn}.\nonumber
\end{eqnarray}
\end{corollary}

\section{Asymptotics for the sample RSPCC and variates}\label{sec4}

The basic ingredients for the asymptotic distribution of the sample
RSPCC and its variates are the weak limits of the associated operators
$\widehat {R}_1$ and $\widehat {R}_2$ (cf. (\ref{e2.17}) and (\ref{e2.18}))
from which these quantities are derived. These limits follow
rather routinely with the help of Corollary \ref{cor3.1}.
It has already been observed that $R_j, \widehat {R}_j\in\mathcal
{L}(\mbox{HS})$ for $j=1,2$.

Let us introduce the following zero-mean Gaussian elements of $\mathcal
{L}(\mbox{HS})$:
%
\begin{eqnarray}
\EuScript{R}_{11} &=& (\varphi '_{1,1}\mathcal{G}_{11}) \Sigma
_{12}\varphi _2(\Sigma _{22})\Sigma _{21}\varphi _1(\Sigma _{11}),\label{e4.1}\\
\EuScript{R}_{12} &=& \varphi _{1}(\Sigma _{11}) \mathcal{G}_{12}
\varphi _2(\Sigma _{22})\Sigma _{21}\varphi _1(\Sigma _{11}),\label{e4.2}\\
\EuScript{R}_{13} &=& \varphi _{1}(\Sigma _{11}) \Sigma _{12}
(\varphi _{2,2}'\mathcal{G}_{22})\Sigma _{21}\varphi _1(\Sigma _{11}),\label{e4.3}\\
\EuScript{R}_{14} &=& \varphi _{1}(\Sigma _{11}) \Sigma _{12}
\varphi _2(\Sigma _{22}) \mathcal{G}_{21}\varphi _1(\Sigma _{11}),\label{e4.4}\\
\EuScript{R}_{15} &=& \varphi _{1}(\Sigma _{11}) \Sigma _{12}
\varphi _2(\Sigma _{22})\Sigma _{21} \varphi
_{1,1}'(\mathcal{G}_{11}),\label{e4.5}\\
\EuScript{R}_1&=&\sum_{j=1}^5 \EuScript{R}_{1j}\label{e4.6}
\end{eqnarray}
and, similarly,
%
\begin{eqnarray}
\EuScript{R}_{21} &=& (\varphi '_{1,2}\mathcal{G}_{22}) \Sigma
_{21}\varphi _2(\Sigma _{11})\Sigma _{12}\varphi _1(\Sigma _{22}),\label{e4.7}\\
\EuScript{R}_{22} &=& \varphi _{1}(\Sigma _{22}) \mathcal{G}_{21}
\varphi _2(\Sigma _{11})\Sigma _{12}\varphi _1(\Sigma
_{22}),\label{e4.8}\\
\EuScript{R}_{23} &=& \varphi _{1}(\Sigma _{22}) \Sigma _{21}
(\varphi _{2,1}'\mathcal{G}_{11})\Sigma _{12}\varphi _1(\Sigma _{22}),\label{e4.9}\\
\EuScript{R}_{24} &=& \varphi _{1}(\Sigma _{22}) \Sigma _{21}
\varphi _2(\Sigma _{11}) \mathcal{G}_{12}\varphi _1(\Sigma _{22}),\label{e4.10}
\\
\EuScript{R}_{25} &=& \varphi _{1}(\Sigma _{22}) \Sigma _{21}
\varphi _2(\Sigma _{11})\Sigma _{12} \varphi
_{1,2}'(\mathcal{G}_{22}),\label{e4.11}\\
\EuScript{R}_2&=&\sum_{j=1}^5 \EuScript{R}_{2j}.\label{e4.12}
\end{eqnarray}
%
%
\begin{theorem}\label{thm4.1}
Let (\ref{e2.1}) be satisfied. We have
%
\begin{equation}
\label{e4.13}
\sqrt{n}(\widehat {R}_j -R_j
)\mathop{\stackrel{d}{\rightarrow}}\EuScript
{R}_j,\qquad\mbox{as }  n\rightarrow \infty , \mbox{ in }
\mathcal{L}(\mathrm{HS})  \mbox{ for } \ j=1,2.
\end{equation}
\end{theorem}
\begin{pf*}{Proof}
It suffices to prove (\ref{e4.13}) for $j=1$. The left-hand side of
(\ref{e4.13}) can be decomposed as $\sum_{j=1}^5 \widehat
{\EuScript{R}}_{1j}$, where, for instance,
%
\begin{equation}
\label{e4.14}
\widehat {\EuScript{R}}_{11} = \sqrt{n} \{ \varphi
_1(\widehat {\Sigma }_{11}) - \varphi _1(\Sigma _{11}) \}
\widehat {\Sigma }_{12}\varphi _2(\widehat {\Sigma
}_{22})\widehat {\Sigma }_{21}\varphi _1(\widehat {\Sigma }_{11}).
\end{equation}
It follows from (\ref{e3.14}) that the first factor in (\ref{e4.14})
equals $\varphi '_{1,2} \mathcal{G}_{22}+\mathscr{o}_p(1)$. Relation
(\ref{e3.3}) and the continuity of the functions in (\ref{e3.14})
imply that the product of the remaining four factors equals $\Sigma
_{21} \varphi _2(\Sigma _{11})\Sigma _{12} \varphi _{1}(\Sigma
_{22}) +\mathscr{O}_p(1)$. In combination, these results yield that
$\widehat {\EuScript{R}}_{11} =\EuScript{R}_{11}+\mathscr{O}_p(1)$.
In a similar manner, one can deal with $\widehat {\EuScript{R}}_{12}, \
\ldots, \widehat {\EuScript{R}}_{15}$.
Eventually, this produces $\sqrt{n}(\widehat {R}_1 - R_1) = \sum
_{j=1}^5 \widehat {\EuScript{R}}_{1j} + \mathscr{o}_p(1)$ and we are done.
\end{pf*}

To establish (\ref{e4.13}), we have exploited the delta-method of
(\ref{e3.14}), based on the Fr\'echet\hspace*{2pt} derivative, in order to deal
with the factors in the product defining $\widehat {R}_j$. Once the
limiting distributions of the random operators have been established,
we may proceed as in Dauxois \textit{et al.} (\citeyear{daporo82}) to find the asymptotic
distributions of eigenvalues and eigenvectors.
For completeness, the required perturbation results in the infinite-dimensional
situation are briefly summarized in the \hyperref[secA]{Appendix} and
proofs of the two main theorems below are included.

First, some more notation will be needed. The compact operators $R_j$
and $\widehat {R}_j$ are nonnegative Hermitian and have spectral
representations
%
\begin{equation}
\label{e4.15}
R_j=\sum_{k=1}^\infty \rho _{jk} Q_{jk},\qquad \widehat {R}_{j} =
\sum_{k=1}^\infty \widehat {\rho }_{jk} \widehat {Q}_{jk},\qquad j=1,2,
\end{equation}
where
$\rho _{j1}>\rho _{j2}> \cdots\downarrow0$ and $\widehat {\rho
}_{j1}>\widehat {\rho }_{j2}> \cdots\downarrow0$ are the distinct
eigenvalues and
$Q_{jk}$, $\widehat {Q}_{jk}$ the orthogonal projections onto the
corresponding finite dimensional eigenspaces.
Assumption (\ref{e2.20}) implies that
%
\begin{equation}
\label{e4.16}
\rho _{j1} = \rho ^2,\qquad Q_{j1} = f^*_j\otimes f^*_j,\qquad j=1,2.
\end{equation}
We also have, by Definition \ref{def2.1} and Theorem \ref{thm2.1}, that
%
\begin{equation}
\label{e4.17}
\widehat {\rho }_{j1} = \widehat {\rho }^2,\qquad j=1,2.
\end{equation}
The operators
%
\begin{equation}
\label{e4.18}
A_j = \sum_{k=2}^\infty \frac{\rho _{j1}}{\rho _{j1}-\rho
_{jk}}Q_{jk},\qquad j=1,2,
\end{equation}
will also be needed.
\begin{theorem}\label{thm4.2}
Let (\ref{e2.1}) and (\ref{e2.19}) be satisfied. The sample
RSPCC then has a normal distribution in the limit:
%
\begin{equation}
\label{e4.19}
\sqrt{n} (\widehat {\rho }^2-\rho ^2 ) \mathop{\stackrel{d}{\rightarrow}}
N(0,\sigma ^2)\qquad \mbox{as } \ n\rightarrow \infty ,
\end{equation}
where
%
\begin{equation}
\label{e4.20}
\sigma ^2 = \mathbb{E}\ \langle {\EuScript{R}_j f^*_j},{f^*_j}
\rangle ^2,\qquad j=1,2.
\end{equation}
\end{theorem}
\begin{pf*}{Proof}
The proof is in the same vein as that of Theorem \ref{thm3.1}.
However, let us now consider the random perturbation $\widehat
{\EuScript{P}} = \widehat {R}_j-R_j$ and define $\Omega _n$ for the
same $\epsilon _n$, but with $\widehat {\EuScript{P}}$ as above. In
the present situation, it is (\ref{e4.13}) that guarantees that
$\mathbb{P}(\Omega _n)\rightarrow 1$ as $n\rightarrow \infty $.

It follows from Theorem \ref{thmA.2} that
%
\begin{equation}
\label{e4.21}
\widehat {Q}_{j1} \mathbf{1}_{\Omega _n} = \widehat {f}^*_j\otimes
\widehat {f}^*_j\mathbf{1}_{\Omega _n}
\end{equation}
for $n$ sufficiently large. Application of Theorem \ref{thmA.3} yields
\begin{eqnarray}
\label{e4.22}
\sqrt{n} (\widehat {\rho }_{j1}-\rho _{j1} )& =& \sqrt
{n} (\widehat {\rho }_{j1}-\rho _{j1} )\mathbf
{1}_{\Omega _n} +\sqrt{n} (\widehat {\rho }-\rho
)\mathbf{1}_{\Omega _n^c}\nonumber \\
& =& \sqrt{n} \langle {\widehat {\EuScript{P}}f^*_j},{f^*_j} \rangle
\mathbf{1}_{\Omega _n} + \EuScript{O}(\|\widehat {\EuScript{P}}\|
^2\mathbf{1}_{\Omega _n}) + \mathscr{o}_p(1)\nonumber\\[-8pt]\\[-8pt]
& =& \langle {\sqrt{n}(\widehat {R}_j-R_j)f^*_j},{f^*_j} \rangle +
\mathscr{o}_p(1)\nonumber \\
&\mathop{\stackrel{d}{\rightarrow}}& \langle {\EuScript{R}_j f^*_j},{f^*_j}
\rangle,\qquad \mbox{as } \ n\rightarrow \infty . \nonumber
\end{eqnarray}
Because of (\ref{e4.16}) and (\ref{e4.17}), the expression on the left
in (\ref{e4.22}) equals the one on the left in (\ref{e4.19}), so
the theorem follows.
\end{pf*}
\begin{theorem}\label{thm4.3}
Assuming the validity of (\ref{e2.1}) and (\ref{e2.19}), we have
%
\begin{equation}
\label{e4.23}
\sqrt{n}(\widehat {f}^*_j -f^*_j)\mathop{\stackrel{d}{\rightarrow}}
A_j\EuScript{R}_jf^*_j,\qquad \mbox{as } \ n\rightarrow \infty , \mbox{ in }  \mathbb{H} \mbox{ for }  j=1,2.
\end{equation}
\end{theorem}
\begin{pf*}{Proof}
Let us consider the same random perturbation $\widehat {\EuScript{P}}
= \widehat {R}_{j} -R_{j}$ and the same sets $\Omega _n$ as in the
proof of Theorem \ref{thm4.2}. Let us also recall (\ref{e4.21}). It
follows from Theorem \ref{thmA.2} that
%
\begin{equation}
\label{e4.24}
\widehat {f}^*_j \mathbf{1}_{\Omega _n} = (f^*_j +A_j \EuScript
{P}f^*_j) \mathbf{1}_{\Omega _n} + \EuScript{O}(\EuScript{P}^2
\mathbf{1}_{\Omega _n}).
\end{equation}
In the same manner as (\ref{e4.22}), we now obtain
%
\begin{equation}
\label{e4.25}
\sqrt{n}(\widehat {f}^*_j -f^*_j) = A_j\sqrt{n} (\widehat
{R}_j-R_j)f^*_j +\mathscr{o}_p(1) \mathop{\stackrel{d}{\rightarrow}}
A_j\EuScript{R}_jf^*_j\qquad \mbox{as } \ n\rightarrow \infty ,
\end{equation}
which proves the theorem.
\end{pf*}

%

\section{Further specification of limiting distributions}\label{sec5}

The distributions on the right in (\ref{e4.19}) and (\ref{e4.23})
contain unknown parameters that must be estimated for practical
implementation. Let us first consider the variance in (\ref{e4.20}).
Substituting (\ref{e4.6}) or (\ref{e4.7}) yields
%
\begin{equation}
\label{e5.1}
\sigma ^2 = \sum_{k=1}^5\sum_{m=1}^5 \mathbb{E} \langle {\EuScript
{R}_{jk}f^*_j},{f^*_j} \rangle \langle {\EuScript
{R}_{jm}f^*_j},{f^*_j} \rangle .
\end{equation}
Subsequent substitution of the expressions for the $\EuScript{R}_{jk}$
shows, after reworking the inner products, that the expression for
$\sigma ^2$ in (\ref{e5.1}) is a sum of terms of the type
%
\begin{equation}
\label{e5.2}
\mathbb{E} \langle {\mathcal{G}f},{g} \rangle \langle {\mathcal
{G}p},{q} \rangle ,
\end{equation}
where $f$, $g$, $p$, $q\in\mathbb{H}$ depend on $\Sigma $ and where
$\mathcal{G}$ is given in (\ref{e3.1}).
\begin{lemma}\label{lem5.1}
If $f$, $g$, $p$, $q\in\mathbb{H}$ are known, we can express (\ref{e5.2}) as
%
\begin{equation}
\label{e5.3}
\mathbb{E} \langle {\mathcal{G}f},{g} \rangle \langle {\mathcal
{G}p},{q} \rangle =\langle {q\otimes p},{\Sigma _{\mathrm{HS}}
g\otimes f} \rangle _{\mathrm{HS}},
\end{equation}
where $\Sigma _{\mathrm{HS}}$ is the covariance operator of $\mathcal
{G}$ in (\ref{e3.2}).
\end{lemma}
\begin{pf*}{Proof}
Let us assume that $f$, $g$, $p$, $q\neq0$ because, otherwise, (\ref{e5.3})
 is trivial. Hence, we can construct two orthonormal bases of
$\mathbb{H}$, viz. $e_1, e_2, \ldots$ and $d_1,d_2, \ldots,$ with
%
\begin{equation}
\label{e5.4}
e_1 = \frac{f}{\|f\|},\qquad d_1 = \frac{p}{\|p\|}.
\end{equation}
Rewriting and evaluating the right-hand side of (\ref{e5.3}), we obtain
\begin{eqnarray}
\label{e5.5}
\langle {q\otimes p},{\Sigma _{\mathrm{HS}} g\otimes f} \rangle
_{\mathrm{HS}} & =& \mathbb{E}\  \langle {q\otimes p},{(\mathcal
{G}\otimes_{\mathrm{HS}} \mathcal{G}) g\otimes f} \rangle _{\mathrm
{HS}}\nonumber \\
& =& \mathbb{E}  \langle {\mathcal{G}},{g\otimes f} \rangle
_{\mathrm{HS}} \langle {\mathcal{G}},{q\otimes p} \rangle _{\mathrm
{HS}} \nonumber\\[-8pt]\\[-8pt]
&=&\mathbb{E}\Biggl\{ \sum_{k=1}^\infty \langle {f},{e_k} \rangle
\langle {\mathcal{G}e_k},{g} \rangle \Biggr\} \Biggl\{ \sum
_{m=1}^\infty \langle {p},{d_m} \rangle \langle {\mathcal
{G}d_m},{q} \rangle \Biggr\} \nonumber \\
&=& \mathbb{E} \langle {\mathcal{G}f},{g} \rangle \langle
{\mathcal{G}p},{q} \rangle ,\nonumber
\end{eqnarray}
as was to be shown.
\end{pf*}

Since the $f$, $g$, $p$, $q$ depend on $\Sigma $, we can replace them
on the right in (\ref{e5.3}) with estimators obtained by substituting
$\widehat {\Sigma }$ for $\Sigma $. Also, $\Sigma _{\mathrm{HS}}$
is unknown and we may replace this operator with the estimator
%
\begin{equation}
\label{e5.6}
\widehat {\Sigma }_{\mathrm{HS}} = \frac{1}{n} \sum_{i=1}^n
[ \{ (X_i-\overline{X})\otimes(X_i-\overline{X}) -\widehat
{\Sigma }\}\otimes_{\mathrm{HS}}\{ (X_i-\overline
{X})\otimes(X_i-\overline{X}) -\widehat {\Sigma }\} ].
\end{equation}

Let us next turn to the Gaussian random element in $\mathbb{H}$, on
the right in (\ref{e4.23}). Substitution of (\ref{e4.6}) or (\ref{e4.7})
 shows that the covariance operator of this random element is
determined by covariances of the type
%
\begin{equation}
\label{e5.7}
\sigma ^2(f,g) = \sum_{k=1}^5\sum_{m=1}^5 \mathbb{E} \langle
{f},{\EuScript{R}_{jk}f^*_j} \rangle \langle {\EuScript
{R}_{jm}f^*_j},{g} \rangle
\end{equation}
and this can be seen to be a sum of terms of type
%
\begin{equation}
\label{e5.8}
\mathbb{E} \langle {p},{\mathcal{G}f^*_j} \rangle \langle
{\mathcal{G}f^*_j},{q} \rangle ,
\end{equation}
in the same way as above.
In this case, explicit expressions for $p$ and $q$ involve the operator
$A_j$ and hence the unknown $\rho _{jk}$ and $Q_{jk}$ (see (\ref{e4.18})
 and (\ref{e4.15})). These quantities can be estimated by the
corresponding quantities for $\widehat {\EuScript{R}}_j$ and $\Sigma
_{\mathrm{HS}}$ can again be estimated by (\ref{e5.6}) so that, in
principle, an estimator of (\ref{e5.7}) is available.
An alternative to this estimation scheme could perhaps be formulated
using resampling and bootstrap methods. We will not explore this idea
further here.

\section{Example and some remarks}\label{sec6}

The purpose of this paper is to establish some fundamental results
regarding functional canonical correlations and their variates, at a
fixed, but arbitrary, level of the regularization parameter $\alpha >0$.
Although the question of how to choose this parameter in practice is
certainly of great interest and relevance, it is not the main
concern of this paper and would require a lengthy discussion of further theory
and numerical simulations beyond the scope and purpose of this work.

As a compromise, in this section, we present an explicit example
that seems suitable for such simulations. It concerns two dependent
standard Brownian motion processes that allow for canonical
correlations in the entire range from $0$ to $1$. To construct these
processes, let
%
\begin{eqnarray}
e_m(t) &=& \sqrt{2} \sin\biggl( \biggl(m-\frac{1}{2}\biggr)\curpi
t\biggr),\qquad t\in\mathbb{R},\ \ m\in\mathbb{N},\label{6.1}\\
\lambda _m &= &\biggl\{ \frac{1}{(m-1/2)\curpi}\biggr\}^2,\qquad\hspace*{23pt} m\in\mathbb{N}.\label{6.2}
\end{eqnarray}
Let $\xi_{jm}$ be i.i.d. $N(0,1)$-random variables for $m\in
\mathbb{N}$ and $j=1,2$. Choose $a_m$, $b_m\in\mathbb{R}$ such that
%
\begin{equation}
\label{6.3}
a_m^2+b_m^2=1,\qquad m\in\mathbb{N},
\end{equation}
and define ($j=1,2$)
%
\begin{eqnarray}
&&\hspace*{10pt} e_{jm} = e_m\bigl(t-(j-1)\bigr) \mathbf{1}_{[j-1,j]} (t) ,\qquad 0\leq
t,\label{6.4}\\
&&\cases{U_{1m} = \sqrt{\lambda _m} \xi_{1m}, &\quad\hspace*{16pt}$ m\in\mathbb
{N}$,\cr
U_{2m} = \sqrt{\lambda _m} ( a_m \xi_{1m} +b_m \xi_{2m}
),&\quad\hspace*{16pt}$ m\in\mathbb{N}$.}\label{6.5}
\end{eqnarray}
For both values of $j$, the
%
\begin{equation}
\label{6.6}
U_{jm} \mbox{ are independent } N(0,\lambda _m),\qquad m\in
\mathbb{N}.
\end{equation}

Obviously,
%
\begin{equation}
\label{6.7}
X_j(t) = \sum_{m=1}^\infty U_{jm} e_{jm}(t), \qquad 0\leq t\leq2,
\end{equation}
is the Karhunen--Lo\`{e}ve expansion of a standard Brownian motion,
starting at $t=0$ for $j=1$ and at $t=1$ for $j=2$. If we define
%
\begin{equation}
\label{6.8}
X(t) = X_1(t)+X_2(t),\qquad 0\leq t\leq2,
\end{equation}
then this process is a random function in $\mathbb{H}=L^2(0,2)$ and $X_j$
can be considered as its projection onto $\mathbb{H}_j = L^2(j-1,j)$.

Because
%
\begin{equation}
\label{6.9}
\gamma _{km} = \mathbb{E}U_{1k} U_{2m} = \sqrt{\lambda _k\lambda
_m}\, a_m \delta _{km},
\end{equation}
a straightforward, but tedious, calculation (see \citeyear{cueugiru06})
shows that $\rho ^2$ is the largest eigenvalue of the diagonal
matrix $\EuScript{R}$ with elements
%
\begin{equation}
\label{6.10}
\EuScript{R}(k,j) =
\cases{\dfrac{a_k^2 \lambda _k^2}{(\alpha +\lambda _k)^2}, &\quad
$k=j$, \cr
0, &\quad $k\neq j$.}
\end{equation}
If we assume that
%
\begin{equation}
\label{6.11}
1\geq a_1^2 \geq a_2^2 \geq\cdots,
\end{equation}
then the largest eigenvalue of this matrix equals
%
\begin{equation}
\label{6.12}
\rho ^2 = \dfrac{a_1^2 \lambda _1^2}{(\alpha +\lambda _1)^2}.
\end{equation}
Choosing $a_1^2=0$ yields $X_1\perp\!\!\!\!\perp X_2$ and $\rho ^2=0$,
and choosing $a_1^2$ close to $1$ and $\alpha $ close to $0$ yields a
$\rho ^2$ close to $1$.

A sample of size $n$ of processes can be obtained by generating $n$
independent, suitably truncated sets of i.i.d. $N(0,1)$-random
variables and $\widehat {R}_1$ can, in principle, be numerically
approximated, by first approximating $\overline{X}$ and $\widehat
{\Sigma }$ in
(\ref{e2.13}). Finally, this should yield a specific value of
$\widehat {\rho }^2$ and hence of $\sqrt{n}(\widehat {\rho
}^2-\rho ^2)$. This sampling process may be repeated $N$ times. Each
of the $N$ runs yields a value of the standardized empirical canonical
correlation and these values could be summarized in a histogram. All
of this might be repeated for several values of the regularization
parameter $\alpha >0$. Numerical procedures are available, but their
implementation is rather involved. Apart from these simulations, some
criterion should be formulated that yields an optimal value of $\alpha
$ in theory, like the mean integrated square error for curve
estimation. The entire issue of gaining insight into the choice of
regularization parameter seems a topic of independent interest.

\renewcommand\thesection{A.{section}}
\setcounter{section}{0}
\setcounter{equation}{0}
\begin{appendix}
\section*{Appendix: Some perturbation theory}\label{secA}

In this appendix, we briefly summarize some results from perturbation
theory. A more general version of these results can be found in a
technical report by Gilliam \textit{et al.} (\citeyear{gihojiru07}. In slightly different
form, Theorem \ref{thmA.2} and Theorem \ref{thmA.3} can be found in
Dauxois \textit{et al.} (\citeyear{daporo82}). Some monographs on perturbation theory for
operators are Kato (\citeyear{ka66}), Rellich (\citeyear{rellich69}) and Chatelin (\citeyear{ch83}). For
matrices, Theorem A.1 can be found in Bhatia (\citeyear{bh07}). It has already been
observed that the delta-method for functions of matrices can be found
in Ruymgaart and Yang (\citeyear{ruya97}).

All operators considered here map the infinite dimensional, separable
Hilbert space $\mathbb{H}$ into itself. As in the main body of the
paper, the inner product and norm in $\mathbb{H}$ will be denoted by
$\langle {\cdot},{\cdot} \rangle $ and $
\|\cdot\|$, respectively, and we will use $\mathcal{L}_H$ to denote
all bounded Hermitian operators on $\mathbb{H}$, with $\mathcal{C}_H$
denoting
the subspace of all compact Hermitian operators and $\mathcal{C}_H^+$
the subset of all strictly positive Hermitian operators. Without
confusion, the operator norm will also be denoted by $\|\cdot\|$.

Let $T\in\mathcal{C}_H^+$ be arbitrary, but fixed. Such an operator
has a spectral representation of the form
%
\begin{equation}
\label{eA.1}
T = \sum_{j=1}^\infty \lambda _j P_j,
\end{equation}
where $\lambda _1>\lambda _2> \cdots\downarrow0$ are the distinct
eigenvalues in decreasing order and $P_1, P_2, \ldots$ are the
projections onto the corresponding finite-dimensional
eigenspaces.

The operator $T$ will be perturbed with a compact Hermitian operator
$\EuScript{P}\in\mathcal{C}_H$. For $r>0$ ,we will write
%
\begin{equation}
\label{eA.2}
\EuScript{O}(\|\EuScript{P}\|^r)
\end{equation}
to indicate any quantity (operator, vector, number) whose norm or
absolute value is of the indicated order as $\|\EuScript{P}\|
\rightarrow 0$.

The perturbed operator $\widetilde {T} =T+\EuScript{P}$ is no longer
strictly positive, but still $\widetilde {T}\in\mathcal{C}_H$. This
operator has spectral representation
%
\begin{equation}
\label{eA.3}
T = \sum_{j=1}^\infty \widetilde {\lambda }_j \widetilde {P}_j,
\end{equation}
where $\widetilde {\lambda }_1, \widetilde {\lambda }_2, \ldots$
are distinct nonzero eigenvalues such that $|\widetilde {\lambda }_1|
\geq|\widetilde {\lambda }_2| \geq \cdots\downarrow0$, and
$\widetilde {P}_1, \widetilde {P}_2, \ldots$ are the projections onto the
corresponding finite-dimensional eigenspaces.

Furthermore, let $\varphi\dvtx D\rightarrow \mathbb{C}$ be
analytic on the open domain $D\subset\mathbb{C}$, where
%
\begin{equation}
\label{eA.4}
D\supset[-\epsilon , \|T\|+\epsilon ]\qquad \mbox{ for some }
\epsilon >0.
\end{equation}
\begin{theorem}\label{thmA.1}
We have
%
\begin{equation}
\label{eA.5}
\varphi (\widetilde {T}) = \varphi (T) + \varphi _T'\EuScript{P}+
\EuScript{O}(\|\EuScript{P}\|^2),
\end{equation}
where $\varphi _T'\dvtx \mathcal{C}_H\rightarrow \mathcal{L}_H$
is bounded and given by
%
\begin{equation}
\label{eA.6}
\varphi '_T \EuScript{P} = \sum_{j\geq1} \varphi '(\lambda _j)
P_j\EuScript{P}P_j+\mathop{\sum\sum}_{j\neq k} \frac{\varphi
(\lambda _k)-\varphi (\lambda _j)}{\lambda _k-\lambda _j}
P_j\EuScript{P}P_k .
\end{equation}
\end{theorem}
\begin{remark}\label{remA.1}
The double sum in (\ref{eA.6}) is a correction term that is needed
because the increment $\Pi\in \mathcal{C}_H$ is arbitrary and
therefore does not, in general, commute with $T$. This generality is
needed for statistical application, as in Theorem~\ref{thm3.1}; see
also Remark~3.1. If $T$ and $\Pi$ do commute, however, then the double sum
would disappear and we would obtain the much simpler expression
%
\begin{equation}
\label{eA.7}
\varphi '_T \EuScript{P} = \sum_{j\geq1} \varphi '(\lambda _j)
P_j\EuScript{P}=(\varphi '(T))\EuScript{P}.
\end{equation}
In other words, in this case, the Fr\'echet derivative $\varphi '_T$
equals the operator $\varphi '(T)$, obtained by applying the usual
functional calculus with the derivative $\varphi '$ of $\varphi $;
see also Dunford and Schwartz (\citeyear{dusch57}, Theorem~VII.6.10 for commuting operators.
\end{remark}
\begin{theorem}\label{thmA.2}
If the range of $P_1$ is one-dimensional so that $P_1 = p_1\otimes p_1$
for some unit vector $p_1\in\mathbb{H}$, then there exists a unit vector
$\widetilde {p}_1\in\mathbb{H}$ such that $\widetilde {P}_1=
\widetilde {p}_1\otimes\widetilde {p}_1$ for $\EuScript{P}$
sufficiently small. We have, moreover, that
%
\begin{equation}
\label{eA.8}
\widetilde {p}_1 = p_1 +A\EuScript{P}p_1 +\EuScript{O}(\|\EuScript{P}\|^2),
\end{equation}
where $A\dvtx \mathbb{H}\rightarrow \mathbb{H}$ is the bounded operator
%
\begin{equation}
\label{eA.9}
A = \sum_{j=2}^\infty \frac{\varphi (\lambda _1)}{\lambda
_1-\lambda _j}\, P_j.
\end{equation}
\end{theorem}
\begin{theorem}\label{thmA.3}
If the range of $P_1$ is one-dimensional and hence $P_1 = p_1\otimes
p_1$ for some unit vector $p_1\in\mathbb{H}$, then we have
%
\begin{equation}
\label{eA.10}
\widetilde {\lambda }_1 =\lambda _1 +\langle {\EuScript
{P}p_1},{p_1} \rangle +\EuScript{O}(\|\EuScript{P}\|^2).
\end{equation}
\end{theorem}
\end{appendix}

\section*{Acknowledgements}
The authors are grateful to the
editor, and an associate editor and referee for useful comments.
Eubank, Gilliam and Ruymgaart gratefully acknowledge support from NSF
Grant DMS-06-05167; Gilliam was also supported by AFOSR Grant
FA9550-07-1-0214.

\printhistory

\end{document}